\documentclass[11pt]{amsart}
\usepackage{amssymb}
\usepackage{mathrsfs}
\usepackage{amsfonts}

\usepackage{amsthm}
\usepackage{amscd}
\usepackage{amsmath}

\newtheorem{theorem}{Theorem}[section]
\newtheorem{lemma}[theorem]{Lemma}

\newtheorem{conj}[theorem]{Conjecture}

\usepackage{latexsym,amssymb}
\usepackage{amsmath}

\renewcommand{\t}{\, | \,}

 \DeclareMathOperator{\ind}{ind}
 \DeclareMathOperator{\ord}{ord}

 \DeclareMathOperator{\supp}{supp}

\begin{document}

\title{On the Index of Sequences over Cyclic Groups}

\author{Weidong Gao}
\address{Center for Combinatorics, LPMC-TJKLC\\
 Nankai University\\
 Tianjin 300071, P.R. China}
 \email{wdgao1963@yahoo.com.cn}
\author{Yuanlin Li}
\address{Department of Mathematics\\
Brock University \\
St.~Catharines, Ontario \\
Canada L2S 3A1} \email{yli@brocku.ca}
\author{Jiangtao Peng}
\address{Center for Combinatorics, LPMC-TJKLC\\
 Nankai University\\
 Tianjin 300071, P.R. China}
 \email{jtpeng1982@yahoo.com.cn}
 \author{Chris Plyley}
\address{Department of Mathematics\\
Brock University \\
St.~Catharines, Ontario \\
Canada L2S 3A1} \email{cp07rp@brocku.ca}
\author{Guoqing Wang}
\address{Center for Combinatorics, LPMC-TJKLC\\
 Nankai University\\
 Tianjin 300071, P.R. China}
 \email{gqwang1979@yahoo.com.cn}

\keywords{zero-sum sequences, index of sequences, cyclic groups}

\subjclass[2010]{11B30, 11B50, 20K01}

\begin{abstract}
 Let $G$ be a finite cyclic group of order $n \ge 2$. Every sequence $S$ over $G$ can be written in the form $S=(n_1g)\cdot \ldots \cdot (n_lg)$
 where $g\in G$ and $n_1,\ldots, n_l \in [1,\ord(g)]$, and the index $\ind (S)$ of $S$ is defined as the minimum of $(n_1+ \ldots + n_l)/\ord (g)$ over all  $g \in G$ with $\ord (g) = n$.    In this paper we prove that a
sequence $S$ over  $G$ of length $|S| = n$ having an element with multiplicity at least
$\frac{n}{2}$ has a subsequence $T$ with  $\ind (T) = 1$, and if the group order $n$  is
a prime, then the assumption on the multiplicity can be relaxed to
$\frac{n-2}{10}$. On the other hand, if $n=4k+2$ with $k \ge 5$, we
provide an example of a sequence $S$ having length   $|S| > n$ and an element with multiplicity
$\frac{n}{2}-1$ which has no subsequence $T$ with
$\ind (T) = 1$. This disproves a conjecture given twenty years ago by Lemke and
Kleitman.
\end{abstract}

\maketitle

\bigskip
\section{Introduction and Main Results}
\bigskip

Let $G$ be an additively written, finite cyclic group and $g \in G$ with $\ord (g) = |G|$. For a sequence
      \[
      S = (n_1g) \cdot \ldots \cdot (n_lg) \quad \text{over} \ G \,, \quad \text{where} \ l \in
      \mathbb N_0 \quad \text{and} \quad n_1, \ldots, n_l \in [1, n] \,,
      \]
we set
      \[
      \| S \|_g = \frac{n_1+ \ldots + n_l}n  \,,
       \]
and then
\[
      \ind (S) = \min \{ \| S \|_h \mid h \in G \ \text{with} \ \ord (h) = |G|  \} \in \mathbb Q_{\ge 0}
      \]
denotes  the \ {\it index} of $S$. The  index of a sequence is a
crucial invariant in the investigation of (minimal) zero-sum
sequences (resp. of zero-sum free sequences) over cyclic groups. It
was first addressed by Lemke and Kleitman (\cite{Kl-Le89}), used as
key tool by Geroldinger (\cite[page 736]{Ge90d}), and then
investigated by Gao \cite{Ga00c} in a systematical way. Since then
it has found a lot of attention in recent years (see
\cite{C-F-S99,Ch-Sm05a,Ga-Ge09b,Ge09a,L-P-Y-Z10a,Po04a,Sa-Ch07a,
Xi-Yu10a,Yu07a}). We briefly discuss some key results.

If $S$ is a minimal zero-sum sequence, then  $|S|\leq 3$, as well as
$|S|\geq \lfloor\frac{n}{2}\rfloor+2$, implies that $\ind (S)=1$
(see \cite{C-F-S99}, \cite{Sa-Ch07a}, \cite{Yu07a}). In contrast to
that, it was shown that for every $k \in [5,
\lfloor\frac{n}{2}\rfloor+1]$, there is a minimal zero-sum
subsequence $T$ of length $|T|=k$ and with $\ind(T)\geq 2$, and that
the same is true for $k=4$ and $\gcd (n,6)\ne 1$. This lead to the
conjecture that, in case $\gcd (n,6)=1$, every minimal zero-sum
sequence $S$ over $G$ of length $|S|=4$ has $\ind(S)=1$.  Li,
Plyley, Yuan and Zeng \cite{L-P-Y-Z10a} recently proved that this
holds true if $n$ is a prime power, but the general case is still
open.

In 1989, Lemke and Kleitman stated the following conjecture
(\cite[page 344]{Kl-Le89}), which we formulate in the present
language.

\begin{conj}
Let $G$ be a  cyclic group of order $n$, $d$ a  divisor of $n$, and
let $S$ be a sequence over $G$ of length $|S| = n$. Then there
exists a  subsequence $T$ of $S$ and  element $g \in G$ with $\ord (g) = n$ such that
\[
d \, \big| \,  n \| T \|_g \, \big| \, n \,.
\]
In the special case $d=n$, this is equivalent to the existence of a subsequence $T$ with $\ind (T) = 1$.
\end{conj}

\smallskip
Indeed the above is the third of three interesting conjectures stated by Lemke and Kleitman  in
\cite{Kl-Le89}. Their first
conjecture has turned out to be  true for all finite abelian
groups (see \cite{Ge93a}), and the second one is still open.
In this paper we
demonstrate that the above conjecture fails  in general (see Theorem \ref{counterexample}), but that it holds true
under an additional assumption on the highest multiplicity of an element occurring in the sequence.
Here are the main results of the present paper (for any undefined terminology or notation the reader is referred to the
beginning of Section \ref{2}).

\begin{theorem}\label{counterexample}
Let $G$ be a cyclic group of order $n \ge 2$, where $n=4k+2$ for some  $k\geq 5$, and let $g \in G$ with $\ord (g) = n$. Then the sequence
\[
S = g^{\frac{n}{2}-3} \left(\frac{n}{2} g \right)
\left( \big(\frac{n}{2}+1\big) g \right)^{\frac{n}{2}-1}
\left( \big(\frac{n}{2}+2 \big) g \right)^{\lfloor\frac{n}{4}\rfloor-2}
\]
has no subsequence $T$ with  $\ind (T) = 1$.
\end{theorem}

\begin{theorem} \label{integer}
Let $G$ be a cyclic group of order $n \ge 2$ and $S$ be a sequence over $G$ of length $|S| = n$. If $\mathsf h(S)<4$ or $\mathsf h(S)\geq n/2$, then $S$ has a
subsequence $T$ with $\ind(T)=1$ and length $|T|\leq \mathsf h(S)$.
\end{theorem}

\begin{theorem} \label{prime2} Let $G$ be a cyclic group of prime order $p>24318$  and
$S$ be a sequence over $G$ of length $|S| = p$.
If $\mathsf h(S)\geq \frac{p-2}{10}$,  then $S$ has a
subsequence $T$ with $\ind(T)=1$.
\end{theorem}

\smallskip
In Section \ref{2} we summarize our notations and give the proof of Theorem \ref{counterexample}. In the following two sections we provide the proofs of
Theorem \ref{integer} and of Theorem \ref{prime2}. We end the paper with a further conjecture and some open problems (see Section \ref{5}).

\bigskip
\section{Notations and Proof of Theorem \ref{counterexample}} \label{2}
\bigskip

Let $\mathbb N$ denote the set of positive integers, $\mathbb P \subset \mathbb N$ the set of prime numbers, and for rational numbers $a, b \in \mathbb Q$
we set $[a, b] = \{ x \in \mathbb Z \mid a \le x \le b \}$. Let $G$ be an additively written  abelian group and $G_0 \subset G$ a subset. We fix the notation concerning sequences
over $G_0$ (which is consistent with \cite{Ga-Ge06b} and
\cite{Ge-HK06a}). Let $\mathcal F(G_0)$ be the free abelian monoid with basis $G_0$. The
elements of $\mathcal F(G_0)$ are called \ {\it sequences} \ over $G_0$.
We write sequences $S \in \mathcal F (G_0)$ in the form
\[
S =  g_1 \cdot \ldots \cdot g_l = \prod_{g \in G} g^{\mathsf v_g (S)}\,,
\]
where $l \in \mathbb N_0$, $g_1, \ldots, g_l \in G_0$, $\mathsf v_g (S) \in \mathbb N_0$ and $\mathsf v_g (S) = 0$ for almost all $g \in G_0$.
We call \ $|S|=l$ the {\it length} of $S$,  $\sigma (S) = g_1+ \ldots + g_l$ the {\it sum} of $S$, $\mathsf v_g (S)$  the \ {\it multiplicity} \ of $g$ in
$S$, $\supp (S) = \{g \in G \mid \mathsf v_g (S) > 0 \}$ the {\it support} of $S$, and we denote by
\[
\mathsf h (S) =    \max \{ \mathsf v_g (S) \mid g \in G \} \in [0, |S|]
  \qquad  \text{the \ {\it maximum of the multiplicities} \ of \
$S$}\,.
\]
For every group homomorphism $\varphi \colon G \to H$, we set $\varphi ( S) = \varphi (g_1) \cdot \ldots \cdot \varphi ( g_l) \in \mathcal F (H)$,
and if $\varphi$ is the multiplication by some $m \in \mathbb N$, then we set $m S = \varphi (S)$.
We say that $S$ is a {\it zero-sum sequence} if $\sigma (S) = 0$, and it is called a {\it minimal zero-sum sequence} if $\sigma (S) = 0$ but
$\sum_{i \in I} g_i \ne 0$ for all $\emptyset \ne I \subsetneq [1, l]$. Suppose that $G$ is finite cyclic. Then a simple calculation (see \cite[Lemma 5.1.2]{Ge09a}) shows that
\[
\begin{aligned}
\ind (S) & = \min \{ \| S \|_h \mid h \in G \ \text{with} \
      \supp (S) \subset \langle h \rangle \} \\
      & = \min \{ \| S \|_h \mid h \in G \ \text{with} \ \langle \supp (S) \rangle = \langle h \rangle
      \} \,.
\end{aligned}
\]

\medskip
\begin{proof}[Proof of Theorem \ref{counterexample}]
Assume to the contrary that $S$ has a subsequence $T$ with $\ind (T)
= 1$. Then there exists an element $h \in G$ with $\ord (h) = n$
such that $\| T \|_h = 1$. We set
\[
g = j h \qquad \text{and} \qquad  T = g^{x} \left(\frac{n}{2} g \right)^{y}
\left( \big(\frac{n}{2}+1\big) g \right)^{z}
\left( \big(\frac{n}{2}+2 \big) g \right)^{w}
\]
where $j \in [1, n-1]$ with $\gcd (j,n) = 1$,  $x \in [0,  n/2-3]$, $y \in [0,1]$, $z \in [0, n/2 -1]$ and $w \in [0,  n/4 -2]$.
Then
\begin{equation}\label{sum T equvalent 0 mod n}
n \| T \|_g = (x+z+2w)+\frac{n}{2}(y+z+w)\equiv 0 \pmod n.
\end{equation}

\smallskip
\noindent
{\bf Case 1.} $j< \frac{n}{4}$.

Then
\[
T = (jh)^{x} \left(\frac{n}{2} h \right)^{y}
\left( \big(\frac{n}{2}+j\big) h \right)^{z}
\left( \big(\frac{n}{2}+2j \big) h \right)^{w} \,.
\]
Since $\| T \|_h = 1$, we infer that $y+z+w\leq 1$ which implies
that $n \| T \|_g \leq x+(\frac{n}{2}+2)\leq
\frac{n}{2}-3+\frac{n}{2}+2<n$, a contradiction.

\smallskip
\noindent
{\bf Case 2.} $\frac{n}{4}<j< \frac{n}{2}$.

Then
\[
T = (jh)^{x} \left(\frac{n}{2} h \right)^{y}
\left( \big(\frac{n}{2}+j\big) h \right)^{z}
\left( \big(2j - \frac{n}{2} \big) h \right)^{w} \,.
\]
Since $\| T \|_h = 1$, we infer that $x\leq 3$ and $z\leq 1$ which
implies that $n \| T \|_g \leq x+z+2w\leq 3+1+2\,(\lfloor
\frac{n}{4}\rfloor-2)<\frac{n}{2}$. Since $x+z+2w>0$ and again by
$\| T \|_h = 1$, we derive that $x+z+2w \equiv 0 \pmod
{\frac{n}{2}}$, a contradiction.

\smallskip
\noindent
{\bf Case 3.} $\frac{n}{2}<j< \frac{3n}{4}$.

Then
\[
T = (jh)^{x} \left(\frac{n}{2} h \right)^{y}
\left( \big(j - \frac{n}{2}\big) h \right)^{z}
\left( \big(2j - \frac{n}{2} \big) h \right)^{w} \,.
\]
Since $\| T \|_h = 1$, we infer that $x+y+w\leq 1$. We assert that
\begin{equation}\label{x+y+w=1}x+y+w=1.\end{equation} Otherwise, $x=y=w=0$ and $n \| T \|_g=z+\frac{n}{2}z \not \equiv 0 \pmod
{\frac{n}{2}}$,  a contradiction to $n \| T \|_g \equiv 0\pmod n$.
Note that $0<x+z+2w<n$. By \eqref{sum T equvalent 0 mod n}, we have
that
\begin{equation}\label{x+z+2w=n/2}
x+z+2w=\frac{n}{2}
\end{equation}
and
\begin{equation}\label{y+z+w equiv 1 mod 2}
y+z+w\equiv 1\pmod 2.
\end{equation}
By \eqref{x+y+w=1} and \eqref{x+z+2w=n/2}, we have $y+z+w\equiv
z+w-y=\frac{n}{2}-1\equiv 0\pmod 2$, a contradiction to \eqref{y+z+w
equiv 1 mod 2}.

\smallskip
\noindent
{\bf Case 4.} $\frac{3n}{4}<j<n$.

Then
\[
T = (jh)^{x} \left(\frac{n}{2} h \right)^{y} \left( \big(j -
\frac{n}{2}\big) h \right)^{z} \left( \big(2j - \frac{3n}{2} \big) h
\right)^{w} \,.
\]
Since $\| T \|_h = 1$, we infer that $x\leq 1$ and $z\leq 3$ which
implies that $n \| T \|_g \leq x+z+2w\leq1+3 + 2\,(\lfloor
\frac{n}{4}\rfloor-2)<\frac{n}{2}$. Clearly, $x+z+2w>0$. From
\eqref{sum T equvalent 0 mod n}, we derive a contradiction.
\end{proof}

\bigskip
\section{Proof of Theorem \ref{integer}}
\bigskip

We need the following two  results. A simple proof of the first one can be found in \cite[Proposition 4.2.6]{Ge09a}
(for historical
comments see \cite{Ha08c}), and a proof of Lemma \ref{Lemma of Ponomarenko} is given in \cite{Po04a}.

\begin{lemma}\label{Lemma of Alon}
Let $G$ be a finite cyclic group and  $S$ be a sequence over $G$ of length $|S| \ge |G|$. Then S
has a zero-sum subsequence $T$ of length $|T|\in [1,\mathsf h(S)]$.
\end{lemma}

\begin{lemma}\label{Lemma of Ponomarenko}
Let $G$ be a finite cyclic group and
$S$ be a minimal zero-sum sequence over $G$ of length $|S|\in [1,3]$.
Then $\ind(S)=1$.
\end{lemma}

\medskip
\begin{proof}[Proof of Theorem \ref{integer}]
We set $n = |G|$ and $h = \mathsf h (S)$. If $h < 4$, then the
assertion follows from Lemmas 3.1 and 3.2. Suppose that $h \ge n/2$.
Let $g \in G$ with $\mathsf v_g (S) = h$. If $\ord (g) < n$, then
$\ord (g) \le n/2 \le h$, and $T = g^{\ord (g)}$ has the required
properties. If $0 \mid S$, then $T = 0$ has the required properties.

Suppose that $\ord (g) = n$ and that $0 \nmid S$. Then we can write $S$ in the form
\[
S = g^h (b_1g) \cdot \ldots \cdot (b_{n-h}g) \quad \text{where} \quad b_1, \ldots, b_{n-h} \in [2, n-1] \,.
\]
Assume to the contrary that $S$ has no subsequence $T$ with the required properties. We continue with the following assertion.

\begin{enumerate}
\item[{\bf A.}] For  every subset $I \subset [1, n-h]$ we have $\sum_{i \in I} b_i \le n-h + |I| - 1 \,.$
\end{enumerate}
If {\bf A} holds, then we apply it with $I = [1, n-h]$ and obtain that
\[
\sum_{i=1}^{n-h} b_i \le 2(n-h)-1 \,,
\]
a contradiction to $b_1, \ldots, b_{n-h} \in [2, n-1]$. We prove {\bf A} by induction on $|I|$. If there were an $i \in [1, n-h]$ such that $b_i \ge n-h+1$,
then $T = g^{n-b_i}(b_i g)$ were a subsequence of $S$ with $\ind (T) = 1$ and length $|T| = n-b_i+1 \le h$, a contradiction. Let $I \subset [1, n-h]$ with $|I| = k+1 \ge 2$,
say $I = [1,k+1]$, and suppose that {\bf A} holds for all proper subsets of $I$. We set $\beta = b_1+ \ldots + b_{k+1}$. By induction hypothesis we get
$\beta - b_i \le n-h+k-1$ for every $i \in [1, k+1]$, which implies that
\[
\beta = \frac{1}{k}(k \beta) = \frac{1}{k} \sum_{i=1}^{k+1} (\beta - b_i) \le \frac{(k+1)(n-h+k-1)}{k} \le n \,
\]
(to get the last inequality, use that $h \ge n/2$ and $k \le n-h-1$).
Thus, if $\beta \ge n-h+k+1$, then $T = g^{n-\beta}(b_1g) \cdot \ldots \cdot (b_{k+1}g)$ is a subsequence of $S$ with $\ind (T) = 1$ and length $|T| = n-\beta + k + 1 \le h$. This is a contradiction, and thus {\bf A} is proved.
\end{proof}

\medskip
Note that the sequence $S$ given in Theorem \ref{counterexample} satisfies
$\mathsf h (S) = \frac{n}{2}-1$. Thus the assumption in Theorem \ref{integer}, that $\mathsf h(S)\geq \frac{n}{2}$,  cannot be weakened for $n\equiv 2 \pmod 4$.

\bigskip
\section{ Proof of Theorem \ref{prime2}}
\bigskip

We fix
our notations which remain valid throughout the whole section.
Let $G$ be a prime cyclic group of order $|G| = p > 24318$, $G^{\bullet} = G \setminus \{0\}$, and let $S$ be a sequence over $G^{\bullet}$ of length $|S| = p$.
If $g \in G^{\bullet}$, $A \subset \mathbb Z$ and $S = (n_1g) \cdot \ldots \cdot (n_lg)$ with $n_1, \ldots, n_l \in [1,p-1]$, then we set
\[
S (A,g) = \prod_{i \in [1,l], n_i \in A} (n_ig) \,.
\]
For an element $g \in G^{\bullet}$, we set
\[
\Sigma_g (S) = \{ p \, \| T \|_g \ \mid  \ T \ \text{is a subsequence of} \ S \ \text{with} \ \| T \|_g \le 1 \} \,,
\]
and we denote by $\mathsf m_g (S)$ the maximal $t \in [1, p]$ such that $\Sigma_g (T) = [1, t]$ for some subsequence $T$ of $S$.
We define
\[
\mathsf m (S) = \max \{ \mathsf m_g (S) \mid g \in G^{\bullet} \} \,.
\]

\medskip
\centerline{\it From now on we fix an element  $g \in G^{\bullet}$ such that $\mathsf m_g (S) = \mathsf m (S)$.}
\medskip

\medskip
\begin{lemma}\label{Lemma of x geq M+2}
Let  $T$ be a subsequence of $S$ such that $\Sigma_g (T) = [1, \mathsf m (S)]$.
Then
$|T|\leq \mathsf m (S)$,  and if $x \in [1, p-1]$ such that $(xg) \t ST^{-1}$, then   $x \geq \mathsf m (S) +2$.
Furthermore, if $\mathsf m (S) = p$, or if there exists an $x \in [1, p-1]$ such that   $(xg) \t ST^{-1}$ and
$x \geq p- \mathsf m (S)$, then $S$ has a subsequence with index $1$.
\end{lemma}

\begin{proof}
By definition, we have $|T| \leq p \, \| T \|_g = \mathsf m (S)$. If
there is some  $x \in [1, p-1]$ with  $(xg) \t ST^{-1}$ and $x \leq
\mathsf m (S)+1$, then $\Sigma_g ( (xg)T) =[1,\min\{p,\mathsf m (S)
+ x\}]$, a contradiction to the maximality of $\mathsf m (S)$. The
second part of this lemma is clear.
\end{proof}

\medskip
\centerline{\it From now on we suppose that  $S$ has  no subsequence with index $1$.}
\medskip

Let $k\geq 2$ be a positive integer, and let
$F[\frac{1}{k}, \frac{k-1}{k}]$ be all irreducible fractions between
$\frac{1}{k}$ and  $\frac{k-1}{k}$ and with denominators in $[2,k]$,
i.e.,
$$
F\left [\frac{1}{k}, \frac{k-1}{k} \right ]=\left \{\frac{a}{b} \ \Big| \ a \in \mathbb N, \, b \in [2, k] \ \text{with} \
\gcd (a,b)=1 \ \text{and} \  \frac{1}{k}\leq \frac{a}{b}\leq \frac{k-1}{k}  \right \}.
$$

\medskip
\begin{lemma}\label{Lemma for two adjacet fractions}
Let $\frac{a}{b}$ and $\frac{c}{d}$ be two adjacent fractions in
$F[\frac{1}{k}, \frac{k-1}{k}]$ with $\frac{a}{b}<\frac{c}{d}$. Then we have
\begin{enumerate}
\item $b+d\geq k+1$.

\smallskip
\item $bc-ad=1$.
\end{enumerate}
\end{lemma}

\begin{proof}
1.  Note that $\frac{a}{b}<\frac{a+c}{b+d}<\frac{c}{d}$.
Since $\frac{a}{b}$ and $\frac{c}{d}$ are adjacent, it follows that
the irreducible fraction with value $\frac{a+c}{b+d}$ is not in
$F[\frac{1}{k}, \frac{k-1}{k}]$. This forces that $b+d\geq k+1$.

\smallskip
2. Since $\gcd (a,b)=1$, there are two integers $u$ and $v$ such that
$bu+av=1$. Note that $b(u+ma)+a(v-mb)=1$ holds for any integer $m$.
Let $x=u+ma$ and $y=mb-v$. Then, $bx-ay=1$. By choosing $m$ suitably
we may assume that $y\leq k$ and $y+b\geq k+1$. It follows that
$y\geq k+1-b>0$ and $x>0$. From $bx-ay=1$ we get
$$
\frac{x}{y}-\frac{a}{b}=\frac{1}{by}.
$$
If $y>1$, then $\frac{x}{y}$ is a fraction in $F[\frac{1}{k},
\frac{k-1}{k}]$. So, either $\frac{c}{d}=\frac{x}{y}$ and we are
done, or $\frac{c}{d}<\frac{x}{y}$. For the latter case we have
$\frac{1}{by}=\frac{x}{y}-\frac{a}{b}=(\frac{x}{y}-\frac{c}{d})+(\frac{c}{d}-\frac{a}{b})=
\frac{b(dx-cy)+y(cb-ad)}{byd}\geq \frac{b+y}{byd}$. This implies
that $d\geq b+y\geq k+1$, a contradiction.

Now assume that $y=1$ and we must have $b=k$. It follows from
$bx-ay=1$ that $a=kx-1$. Therefore, $x=1$ and $a=k-1$. So,
$\frac{a}{b}=\frac{k-1}{k}$ is the biggest fraction in
$F[\frac{1}{k}, \frac{k-1}{k}]$, a contradiction.
\end{proof}

\medskip
We set
\[
k= \Big\lfloor\frac{p}{\mathsf m (S)}\Big\rfloor \,, \quad f = \Big| \, F \Big[\frac{1}{k}, \frac{k-1}{k} \Big] \, \Big | \,,
\]
and we arrange all
fractions in $F[\frac{1}{k}, \frac{k-1}{k}]$ increasingly; so let
\[
\frac{a_1}{b_1} < \ldots < \frac{a_f}{b_f}
\]
denote the elements of $F[\frac{1}{k}, \frac{k-1}{k}]$. Furthermore, we set
\[
S_1 = S ([1, \mathsf m (S)], g) \, \quad S_2 = S([\mathsf m (S)+2, \frac{p-1}{b_1}], g)
\]
and, for every $i\in
[1,f]$, we set
\[
S_{2i+1}=S \Big( \Big[\frac{a_ip+1}{b_i},
\frac{a_ip+\mathsf m (S)}{b_i} \Big], g \Big) \quad \text{and} \quad S_{2i+2}=S \Big( \Big[\frac{a_ip+\mathsf m (S)+1}{b_i},
\frac{a_{i+1}p-1}{b_{i+1}} \Big], g \Big) \,.
\]
Furthermore, for every $i\in [2,k]$, we define
\[
R_i = S ( \{x \in [1, p] \mid \text{If} \ x_i \in [1, p] \ \text{with} \ p \t (x_i-ix), \ \text{then} \ x_i \in [1,  \mathsf m (S)] \ \text{and} \ \gcd(x_i,i)=1\}, g) \,.
\]

\medskip
\begin{lemma}\label{Lemma for S_i}
We have $S=\prod_{j=1}^{2f+1}S_j$.
\end{lemma}

\begin{proof}
This is clear by construction.
\end{proof}

\medskip
\begin{lemma}\label{Lemma the length of S}
Suppose that
\[
4\leq \mathsf m (S)\leq \frac{p-3}{2} \quad \text{and} \quad
\max \Big\{\frac{p-\mathsf m (S)-2}{\mathsf m (S)}, \ \frac{p-\mathsf m (S)}{\mathsf m (S)+1} \Big\} \leq k \leq \frac{p+1}{\mathsf m (S)} \,.
\]
\begin{enumerate}
\item  $|S_{2i+2}|\leq b_{i+1}-1$ for every $i\in [0,f-1]$.

\smallskip
\item $p=|S| \leq \mathsf m (S)+\sum_{i=2}^k\sum_{j \in [1, i-1] \ \text{with} \ \gcd (i,j)=1}(i-1)+\sum_{i=2}^k|R_i|$.
\end{enumerate}
\end{lemma}

\begin{proof}
1. Suppose that  $i=0$. Then  $S_2 = S ([\mathsf m (S)+2,
\frac{p-1}{b_1}], g)$ and $b_1=k$. If $|S_2|\geq b_1=k$, then we can
take a $k$-term subsequence $U$ of $S_2$. Note that $p-1\geq p
\|U\|_g \geq k(\mathsf m (S)+2)\geq p-\mathsf m (S)$ and one can
find a subsequence $V$ of $S_1$ such that $UV$ has index $1$, a
contradiction.

Now suppose  that $i\in [1, f-1]$, and assume to the contrary that $|S_{2i+2}|\geq b_{i+1}$. We choose  an
arbitrary $b_{i+1}$-term subsequence $X$ of $S_{2i+2}$, and write $b_iS$ in the form
\[
b_iS = \big( x_1g \big) \cdot \ldots \cdot \big( x_p g \big) \quad \text{with} \quad x_1, \ldots, x_p \in [1, p-1] \,.
\]
It follows from Lemma \ref{Lemma for two adjacet
fractions} that $a_{i+1}b_i-a_ib_{i+1}=1$, and so
$b_i(\frac{a_{i+1}p-1}{b_{i+1}})-a_ip=\frac{p-b_i}{b_{i+1}}$.
Thus for every $\nu \in [1,p]$ with $(x_{\nu}g) \t  S_{2i+2}$, we infer that   $x_{\nu} \in [\mathsf m (S)+1,
\frac{p-b_i}{b_{i+1}}]$ and $x_{\nu} \equiv -a_ip \pmod {b_i}$. Therefore we get, since by
Lemma \ref{Lemma for two adjacet fractions} , $b_i+b_{i+1}\geq k+1$,
\[
p-b_i\geq p \|b_iX\|_{g} \geq b_{i+1}(\mathsf m (S)+1)\geq p-b_i\mathsf m (S)
\]
and
\[
p \|b_iX\|_{g} \equiv -b_{i+1}a_ip=(1-a_{i+1}b_i)p\equiv p
\pmod {b_i} \,.
\]
Therefore there exists a subsequence $Y$ of $S_1$ such
that $p \|b_i(X Y) \|_{g}=p$, a contradiction.

\smallskip
2. For every $\ell \in [2,k]$, we have
$R_{\ell}=\prod_{b_i=\ell}S_{2i+1}$, and hence
\[
S=S_1\prod_{i=0}^{f-1}S_{2i+2}\prod_{\ell=2}^kR_{\ell} \,.
\]
Now 2. follows from 1.
\end{proof}

\medskip
\begin{lemma}\label{Lemma n intergers addition}
Let $\ell \in \mathbb{N}_{\geq 2}$ and $S \in \mathcal F ( \mathbb Z
)$ be a sequence of length $|S| = \ell$. Suppose that every element
from $S$ is co-prime to $\ell$. Then for every $m \in \mathbb Z$
there exists a subsequence $S_m$ such that $\sigma (S_m) \equiv m
\pmod \ell$. Moreover, if $m \notin \ell \mathbb Z$, then we get
$S_m \ne S$.
\end{lemma}

\begin{proof}
Let  $\varphi \colon \mathbb Z \to \mathbb Z / \ell \mathbb Z$ be the canonical epimorphism and $\varphi (S) =
a_1 \cdot \ldots \cdot  a_l$.
We denote by $A = \{a_1,0\}+\ldots+\{a_{\ell-1},0\} \subset  \mathbb Z / \ell \mathbb Z$ the sumset, and by
$H = \text{\rm Stab} (A)$ the stabilizer of $A$. Clearly, it suffices to verify that $A = \mathbb Z / \ell \mathbb Z$.  If $H$ would be a proper subgroup of $\mathbb Z / \ell \mathbb Z$, then Kneser's Theorem would imply that
\[
|A| \geq
\sum\limits_{i=1}^{\ell-1}|\{a_i,0\}+H|-(\ell-2)|H|=(\ell-1)
2|H|-(\ell-2)|H|\geq \ell \,,
\]
whence $A = H = Z / \ell \mathbb Z$. Thus $H = \mathbb Z / \ell \mathbb Z$,
which implies that $A = \mathbb Z / \ell \mathbb Z$, and we are done.
\end{proof}

\medskip
\begin{lemma}\label{Lemma of R_t and R_l}
Let $
t, \ell \in [2,  k-1]$ with $t < \ell$ and $d=\gcd(t,\ell)<t$, and let $u\in[2, \mathsf m (S)]$.
If
\[
\frac{(t-d)p-\ell}{t\ell}\leq \mathsf m (S)\leq\frac{dp}{\ell}-t(u-1) \ , \quad \text{then} \quad
|R_t|=0 \quad \text{ or} \quad |R_\ell|\leq \frac{p-\ell \mathsf m (S)-2\ell
+1}{u}+2\ell-1 \,.
\]
\end{lemma}

\begin{proof}
Suppose that $|R_t|>0$. Let $x \in [1, p-1]$ such that  $(xg) \t R_t$, and let $x_{\ell} \in [1, p-1]$ such that $p \t (\ell x - x_{\ell})$. By
the definition of $R_t$, we get
\[
x_{\ell} \in \bigcup_{i \in [1, t-1] \ \text{with} \ \gcd(i,t)=1} \Big[\frac{\ell
ip+\ell}{t},\frac{\ell ip+\ell \mathsf m (S)}{t} \Big] \,,
\]
and thus,
\[
x_{\ell}  \in\bigcup_{i \in [1, t-1] \ \text{with} \ d \t i} \Big[\frac{ ip+\ell}{t},\frac{ ip+\ell \mathsf m (S)}{t} \Big]
 \ \subset \Big[\frac{ dp+\ell}{t},\frac{ (t-d)p+\ell
\mathsf m (S)}{t} \Big] \subset \Big[p-\ell \mathsf m (S), p-\ell(u-1) \Big] \,.
\]
If $|(\ell R_{\ell})([1,u-1],g)|\geq \ell$, then, by Lemma \ref{Lemma n
intergers addition} and the definition of $R_t$, we may choose a
subsequence $W$ of $R_{\ell}$ of length at most $\ell$ such that
$(\ell W)([1,u-1],g) = \ell W$ and $x_{\ell} + p \| \ell
W \|_g \equiv p\pmod \ell$. Since $p \| \ell
W \|_g \leq \ell(u-1)$,
we have $x_{\ell} + p \| \ell
W \|_g \in [p-\ell \mathsf m (S),p]$. Thus, we
can construct a subsequence of $(x g) WS_1$ of index $1$, a
contradiction. Therefore,
\begin{equation}\label{Lemma 2.6 for u-1}
|(\ell R_{\ell})([1,u-1],g)|   \leq \ell-1.
\end{equation}
If $|R_{\ell}|<\ell$ then we are done. Otherwise, by Lemma
\ref{Lemma n intergers addition}, we get a subsequence $R_0$ of
$R_{\ell}$ with $p \| \ell R_0 \|_p \equiv p\pmod \ell$ and
\begin{equation}\label{Lemma 2.6 relation}
|R_0|\geq |R_{\ell}|-\ell.
\end{equation}

\smallskip
We assert that
\begin{equation}\label{Lemma 2.6 for R_0}
p \| \ell R_0 \|_p \leq p-\ell \mathsf m (S)-\ell.
\end{equation}
Assume to the contrary that
$p \| \ell R_0 \|_p \geq p-\ell \mathsf m (S)$, choose $T$ to be the minimal
subsequence of $R_0$ such that $p \| \ell T \|_g \geq p-\ell \mathsf m (S)$ and
$p \| \ell T \|_g \equiv p \pmod \ell$. If $p \| \ell T \|_g \leq
p$, then we can construct a subsequence of $TS_1$ with index $1$, a
contradiction. Now suppose that $p \| \ell T \|_g >p$. If $y \in [1, p-1]$ such that $(yg) \t  R_{\ell}$ and $y_{\ell} \in [1, p-1]$ such that $p \t (\ell y - y_{\ell})$, then $y_{\ell} \in [1,  \mathsf m (S)]$ and $\gcd(
y_{\ell},\ell)=1$. By Lemma \ref{Lemma n intergers addition}, by
dropping at most $\ell$ terms from $T$, we get a proper subsequence
$\tilde{T}$ such that $p \| \ell \tilde{T} \|_g \geq p-\ell \mathsf m (S)$ and
$p \| \ell \tilde{T} \|_g \equiv p\pmod \ell$, a contradiction to
the minimality of $T$. Therefore, (\ref{Lemma 2.6 for R_0}) holds.

\smallskip
By \eqref{Lemma 2.6 for u-1}, we have that $p \| \ell R_0 \|_g \geq
(\ell-1)+u  \big(|R_0|-\ell+1 \big)$. This together with
\eqref{Lemma 2.6 for R_0} gives that $|R_0|\leq \frac{p-\ell \mathsf
m (S)-2\ell +1}{u}+\ell-1$. Now the lemma follows from \eqref{Lemma
2.6 relation}.
\end{proof}

\medskip
\begin{lemma}\label{Lemma the estimate by the intergers sequence}
Let $t\in [2,k]$, and let
$1=\alpha_1<\alpha_2< \ldots$ denote all positive integers coprime to $t$. If
\[
\mathsf m (S)\leq \frac{p-2t+w \alpha_{u+1}+2}{t+\sum\limits_{i=2}^u \alpha_i} \quad  \text{for some} \quad w,u\in \mathbb{N}_0 \,,
\]
then
\[
|R_t|\leq \frac{p-(t+
\sum\limits_{i=2}^u \alpha_i)\mathsf m (S)-2t+2}{\alpha_{u+1}}+ \delta_u(u-1)\mathsf m (S)+2t+w \quad \text{where} \quad
\delta_u = \begin{cases}
           0 \ \text{ for} \ u=0 \\ 1 \ \text{ for} \  u\geq 1
           \end{cases}
\]
\end{lemma}

\begin{proof}
Assume to the contrary that $|R_t|$ is strictly larger than the above bound. Since
\[
\mathsf m (S)\leq \frac{p-2t+w \alpha_{u+1}+2}{t+\sum\limits_{i=2}^u
\alpha_i} \ , \quad \text{it follows that} \quad |R_t|\geq 2t+1 \,.
\]
By Lemma \ref{Lemma n intergers addition}, there exists a
nonempty subsequence $R_0$ of $R_t$ with
\begin{equation}\label{Lemma 2.7 relation}
p \|t R_0 \|_{g}\equiv p\pmod t \quad \text{and} \quad |R_0|\geq |R_t|-t.
\end{equation}
Similarly to Lemma \ref{Lemma of R_t and R_l}, we can prove that
\begin{equation}\label{Lemma 2.7 for R_0}
p \|t R_0 \|_{g}\leq p-t\mathsf m (S)-t.
\end{equation}
Note that $tR_0$ contains $\alpha_1 g = g$ at most $t-2$ times,
because otherwise we would get
\[
\mathsf m(S)\geq \mathsf m_{g} (t S) \geq t\mathsf m_g (S)+t-1>\mathsf m_g (S) = \mathsf m (S) \,,
\]
a contradiction. Since $\mathsf v_{\alpha_i g}(S)\leq \mathsf h(S)\leq \mathsf m (S)$ for all $i\geq
2$, it follows that
\[
p \|t R_0 \|_{g}\geq   \alpha_1(t-2)+
\Big(\sum_{i=2}^u \alpha_i  \Big)\mathsf m (S) +\alpha_{u+1} \big(|R_0|-(u-1)\mathsf m (S)-(t-2) \big) \,.
\]
By \eqref{Lemma 2.7 for R_0}, we have $|R_0|\leq \frac{p-(t+
\sum\limits_{i=2}^ua_i)\mathsf m
(S)-2t+2}{\alpha_{u+1}}+\delta(u-1)\mathsf m (S)+t-2$. By
\eqref{Lemma 2.7 relation}, we derive a contradiction.
\end{proof}

\bigskip
\begin{proof}[Proof of Theorem \ref{prime2}]
We use all the notations which have been fixed at the beginning of this section. In particular, we assume to the contrary that there
exists a sequence  $S \in \mathcal F (G^{\bullet})$ of length $|S| = p$ which
has no subsequence with index $1$. We have to derive a contradiction.

Clearly, we have $\mathsf h (S) \le \mathsf m (S)\leq p-1$.
Lemma \ref{Lemma of x geq M+2} implies that, for every $x \in [1, p-1]$ with $(xg) \t ST^{-1}$,
we have $\mathsf m (S)+2 \leq x \leq p-\mathsf m (S)-1$. Thus it follows that
\[
\frac{p-2}{10} \le \mathsf h (S) \le \mathsf m (S)\leq \frac{p-3}{2} \ .
\]
We distinguish several cases.

\smallskip
\noindent
{\bf Case 1.} $\frac{p-2}{3}\leq \mathsf m (S)\leq \frac{p-3}{2}$.

With $k=2$ in Lemma \ref{Lemma the length of S}, we have
$$
p\leq \mathsf m (S)+1+|R_2|.
$$
Applying Lemma \ref{Lemma the estimate by the intergers sequence}
with $u=0$ and $w=6$, we infer that
$$|R_2|\leq p-2\mathsf m (S)+8.$$
It follows that $p\leq \mathsf m (S)+1+|R_2|=\mathsf m (S)+1+p-2\mathsf m (S)+8<p,$ a contradiction.

\smallskip
\noindent
{\bf Case 2.} $\frac{p+3}{4}\leq \mathsf m (S)\leq \frac{p-4}{3}$.

With $k=3$ in Lemma \ref{Lemma the length of S}, we have
\[
p\leq \mathsf m (S)+1+2+2+|R_2|+|R_3| \,.
\]
Applying Lemma \ref{Lemma the estimate by the intergers sequence}
with $u=1$ and $w=6$, we infer that
\[
|R_2|\leq \frac{p-2\mathsf m (S)+28}{3} \quad \text{and} \quad
|R_3|\leq \frac{p-3\mathsf m (S)+20}{2} \,.
\]
It follows that
\[
p\leq
\mathsf m (S)+5+\sum_{i=2}^3|R_i|=\mathsf m (S)+5+\frac{p-2\mathsf m (S)+28}{3}+\frac{p-3\mathsf m (S)+20}{2}<p \,,
\]
a contradiction.

\smallskip
\noindent
{\bf Case 3.}  $\frac{p-2}{5}\leq \mathsf m (S)\leq \frac{p+1}{4}$.

With $k=4$ in Lemma \ref{Lemma the length of S}, we have
\[
p\leq \mathsf m (S)+1+2\cdot 2+3\cdot 2+|R_2|+|R_3|+|R_4| \,.
\]
Applying Lemma \ref{Lemma the estimate by the intergers sequence}
with $u=1$ and $w=6$, we infer that
\[
|R_2|\leq \frac{p-2\mathsf m (S)+28}{3} \ , \quad
|R_3|\leq \frac{p-3\mathsf m (S)+20}{2} \quad \text{and} \quad
|R_4|\leq \frac{p-4\mathsf m (S)+36}{3} \,.
\]
It follows that
\[
p\leq \mathsf m (S)+11+\frac{p-2\mathsf m (S)+28}{3}+\frac{p-3\mathsf m (S)+20}{2}+\frac{p-4\mathsf m (S)+36}{3}<
p \,,
\]
a contradiction.

\smallskip
\noindent
{\bf Case 4.} $\frac{p-1}{6}\leq \mathsf m (S)\leq \frac{p-3}{5}$.

With $k=5$ in Lemma \ref{Lemma the length of S}, we have
\[
p\leq \mathsf m (S)+27+\sum_{i=2}^5|R_i| \,.
\]
Applying Lemma \ref{Lemma the estimate by the intergers sequence}
with $u=1$ and $w=6$, we infer that
\[
\begin{aligned}
|R_2| & \leq \frac{p-2\mathsf m (S)+28}{3} \ , \quad
|R_3|\leq \frac{p-3\mathsf m (S)+20}{2} \ ,  \\
|R_4| & \leq \frac{p-4\mathsf m (S)+36}{3} \ , \quad
|R_5|\leq \frac{p-5\mathsf m (S)+24}{2} \,.
\end{aligned}
\]
Applying Lemma \ref{Lemma of R_t and R_l} with $t=2,\ell=3$ and
$u=12$, we obtain that either
\[
|R_2|=0 \quad \text{ or} \quad  |R_3|\leq
\frac{p-3\mathsf m (S)+55}{12} \,,
\]
and therefore
\[
|R_2|+|R_3|\leq
\max\{\frac{p-2\mathsf m (S)+28}{3}+\frac{p-3\mathsf m (S)+55}{12},\frac{p-3\mathsf m (S)+20}{2}\}=\frac{5p-11\mathsf m (S)+167}{12} \,.
\]
Summing up we obtain that
\[
\begin{aligned}
p & \leq
\mathsf m (S)+27+\sum_{i=2}^5|R_i|=\mathsf m (S)+27+(|R_2|+|R_3|)+|R_4|+|R_5| \\ & \leq
\frac{5p-11\mathsf m (S)+167}{12}+\frac{p-4\mathsf m (S)+36}{3}+\frac{p-5\mathsf m (S)+24}{2}+27<p \,,
\end{aligned}
\]
a contradiction.

\smallskip
\noindent
{\bf Case 5.} $\frac{p-5}{7}\leq \mathsf m (S)\leq \frac{p-5}{6}$.

With $k=6$ in Lemma \ref{Lemma the length of S}, we have
\[
p\leq \mathsf m (S)+37+\sum_{i=2}^6|R_i| \,.
\]
Applying Lemma \ref{Lemma the estimate by the intergers sequence}
with $u=2$ and $w=0$, we infer that
\[
|R_2|\leq \frac{p+18}{5} \quad \text{and} \quad
|R_3|\leq\frac{p-\mathsf m (S)+20}{4} \,.
\]
Applying Lemma \ref{Lemma the estimate by the intergers sequence}
with $u=1$ and $w=6$, we infer that
\[
|R_4|\leq \frac{p-4\mathsf m (S)+36}{3} \ , \quad
|R_5|\leq \frac{p-5\mathsf m (S)+24}{2} \ , \quad \text{and} \quad
|R_6|\leq \frac{p-6\mathsf m (S)+80}{5} \,.
\]
Summing up we obtain that
\[
\begin{aligned}
p & \leq
\mathsf m (S)+37+\sum_{i=2}^6|R_i| \\ & =\mathsf m (S)+37+\frac{p+18}{5}+\frac{p-\mathsf m (S)+20}{4}+\frac{p-4\mathsf m (S)+36}{3}+\frac{p-5\mathsf m (S)+24}{2}+\frac{p-6\mathsf m (S)+80}{5} \\ & <p \,,
\end{aligned}
\]
a contradiction.

\smallskip
\noindent
{\bf Case 6.} $\frac{p-2}{8}\leq \mathsf m (S)\leq \frac{p-3}{7}$.

With $k=7$ in Lemma \ref{Lemma the length of S}, we have
\[
p\leq \mathsf m (S)+73+\sum_{i=2}^7|R_i| \,.
\]
Applying Lemma \ref{Lemma the estimate by the intergers sequence}
with $u=2$ and $w=0$, we infer that
\[
|R_2| \leq \frac{p+18}{5} \quad \text{and} \quad
|R_3| \leq\frac{p-\mathsf m (S)+20}{4} \,.
\]
Applying Lemma \ref{Lemma the estimate by the intergers sequence}
with $u=1$ and $w=6$, we infer that
\[
\begin{aligned}
|R_4| & \leq \frac{p-4\mathsf m (S)+36}{3} \ , \quad
|R_5|\leq \frac{p-5\mathsf m (S)+24}{2} \ , \\
|R_6| & \leq \frac{p-6\mathsf m (S)+80}{5} \ , \quad
|R_7|\leq \frac{p-7\mathsf m (S)+28}{2} \,.
\end{aligned}
\]
Applying Lemma \ref{Lemma of R_t and R_l}, with $t=2$, $\ell=5$ and
$u=10$, we infer that
\[
|R_2|+|R_5|
\leq\max\{\frac{p-5\mathsf m (S)+4}{2},\frac{p+18}{5}+\frac{p-5\mathsf m (S)-9}{10}+9\}
=\frac{3p-5\mathsf m (S)+117}{10} \,.
\]
Summing up we obtain that
\[
\begin{aligned}
p & \leq
\mathsf m (S)+73+\sum_{i=2}^7|R_i|=\mathsf m (S)+73+(|R_2|+|R_5|)+|R_3|+|R_4|+|R_6|+|R_7| \\ & \leq
\mathsf m (S)+73+\frac{3p-5\mathsf m (S)+117}{10}+\frac{p-\mathsf m (S)+20}{4}+\frac{p-4\mathsf m (S)+36}{3} \\ & \qquad \qquad +\frac{p-6\mathsf m (S)+80}{5}+\frac{p-7\mathsf m (S)+28}{2}<p \,,
\end{aligned}
\]
a contradiction.

\smallskip
\noindent
{\bf Case 7.} $\frac{p-2}{9}\leq \mathsf m (S)\leq \frac{p-3}{8}$.

With $k=8$ in Lemma \ref{Lemma the length of S}, we have
\[
p\leq \mathsf m (S)+111+\sum_{i=2}^8|R_i| \,.
\]
Applying Lemma \ref{Lemma the estimate by the intergers sequence}
with $u=2$ and $w=0$, we infer that
\[
\begin{aligned}
|R_2| & \leq \frac{p+18}{5} \ , \quad
|R_3|\leq\frac{p-\mathsf m (S)+20}{4} \ , \\
|R_4| & \leq\frac{p-2\mathsf m (S)+34}{5} \ , \quad
|R_5|\leq \frac{p-4\mathsf m (S)+22}{3} \,.
\end{aligned}
\]
Applying Lemma \ref{Lemma the estimate by the intergers sequence}
with $u=1$ and $w=6$, we infer that
\[
|R_6|\leq \frac{p-6\mathsf m (S)+80}{5} \ , \quad
|R_7|\leq \frac{p-7\mathsf m (S)+28}{2}  \quad \text{and} \quad
|R_8|\leq\frac{p-8\mathsf m (S)+52}{3} \,.
\]
Applying Lemma \ref{Lemma of R_t and R_l} with $t=2$, $\ell \in \{5,7\}$ and
$u=20$, we can prove that either
\[
|R_2|=0 \quad \text{ or} \quad  |R_i|\leq
\frac{p-i\mathsf m (S)-2i+1}{20}+2i-1 \quad \text{ for } \quad i \in \{5,7 \} \,,
\]
and therefore
\[
\begin{aligned}
|R_2|+|R_5|+|R_7| & \leq \max \Big\{
\frac{p-4\mathsf m (S)+22}{3}+\frac{p-7\mathsf m (S)+28}{2}, \\ & \qquad \qquad \frac{p-\mathsf m (S)+20}{4}+\frac{p-5\mathsf m (S)-9}{20}+9+\frac{p-7\mathsf m (S)-13}{20}+13 \Big\}
 \\ & =\frac{5p-29\mathsf m (S)+128}{6} \,.
\end{aligned}
\]
Applying Lemma \ref{Lemma of R_t and R_l} with $t=4,\ell=6$ and
$u=10$, we obtain that either
\[
|R_4|=0 \quad \text{ or} \quad  |R_6|\leq \frac{p-6
\mathsf m (S)-11 }{10}+11
\]
and therefore
\[
|R_4|+|R_6|\leq
\max \Big\{\frac{p-2\mathsf m (S)+34}{5}+\frac{p-6 \mathsf m (S)-11
}{10}+11,\frac{p-6\mathsf m (S)+80}{5} \Big\}=\frac{3p-10\mathsf m (S)+167}{10} \,.
\]
Summing up we obtain that
\[
\begin{aligned}
p & \leq
\mathsf m (S)+111+\sum_{i=2}^8|R_i|=\mathsf m (S)+111+(|R_2|+|R_5|+|R_7|)+(|R_4|+|R_6|)+|R_3|+|R_8| \\ & \leq
\mathsf m (S)+111+\frac{5p-29\mathsf m (S)+128}{6}+\frac{3p-10\mathsf m (S)+167}{10}+\frac{p-\mathsf m (S)+20}{4}+\frac{p-8\mathsf m (S)+52}{3} \\ & <p \,,
\end{aligned}
\]
a contradiction.

\medskip
{\bf Case 8.} $\frac{p-2}{10}\leq \mathsf m (S)\leq \frac{p-4}{9}$.

With $k=9$ in Lemma \ref{Lemma the length of S}, we have
\[
p\leq \mathsf m (S)+159+\sum_{i=2}^9|R_i| \,.
\]
Applying Lemma \ref{Lemma the estimate by the intergers sequence}
with $u=2$ and $w=0$, we infer that
\[
\begin{aligned}
|R_2|\leq \frac{p+18}{5} \ , \quad
|R_3|\leq\frac{p-\mathsf m (S)+20}{4} \ , \\
|R_4|\leq\frac{p-2\mathsf m (S)+34}{5} \ , \quad
|R_5|\leq \frac{p-4\mathsf m (S)+22}{3} \,.
\end{aligned}
\]
Applying \ref{Lemma the estimate by the intergers sequence} with
$u=1$ and $w=6$, we infer that
\[
\begin{aligned}
|R_6| & \leq \frac{p-6\mathsf m (S)+80}{5} \ , \quad
|R_7|\leq \frac{p-7\mathsf m (S)+28}{2} \ , \\
|R_8| & \leq\frac{p-8\mathsf m (S)+52}{3} \ , \quad
|R_9|\leq \frac{p-9\mathsf m (S)+32}{2} \,.
\end{aligned}
\]
Applying Lemma \ref{Lemma of R_t and R_l} with $t=2,\ell \in \{5,7\}$ and
$u=10$, we obtain that either
\[
|R_2|=0 \quad \text{ or} \quad  |R_i|\leq \frac{p-i
\mathsf m (S)-2i +1}{10}+2i-1 \quad  \text{for} \quad i \in \{5,7\} \,,
\]
and therefore
\[
\begin{aligned}
|R_2|+|R_5|+|R_7| & \leq
\max \Big\{\frac{p-4\mathsf m (S)+22}{3}+\frac{p-7\mathsf m (S)+28}{2}, \\ & \frac{p+18}{5}+\frac{p-5
\mathsf m (S)-9}{10}+9+\frac{p-7 \mathsf m (S)-13}{10}+13 \Big\}  =\frac{5p-29\mathsf m (S)+128}{6} \,.
\end{aligned}
\]
Applying Lemma \ref{Lemma of R_t and R_l} with $t=3,\ell=8$ and
$u=5$, we obtain that either
\[
|R_3|=0 \quad \text{ or} \quad  |R_8|\leq \frac{p-8 \mathsf m (S)-15
}{8}+15 \,,
\]
and therefore
\[
|R_3|+|R_8|\leq
\max \Big\{\frac{p-\mathsf m (S)+20}{4}+\frac{p-8 \mathsf m (S)-15
}{8}+15,\frac{p-8\mathsf m (S)+52}{3} \Big\}=\frac{3p-10\mathsf m (S)}{8}+20 \,.
\]
Summing up we obtain that
\[
\begin{aligned}
p & \leq
\mathsf m (S)+159+\sum_{i=2}^9|R_i| =M+159+(|R_2|+|R_5|+|R_7|)+(|R_3|+|R_8|)+|R_4|+|R_6|+|R_9| \\ & \leq
\mathsf m (S)+159+\frac{5p-29\mathsf m (S)+128}{6}+(\frac{3p-10\mathsf m (S)}{8}+20) \\ & \qquad \qquad + \frac{p-2\mathsf m (S)+34}{5}+\frac{p-6\mathsf m (S)+80}{5}+\frac{p-9\mathsf m (S)+32}{2}<p \,,
\end{aligned}
\]
a contradiction.
\end{proof}

\bigskip
\section{A Conjecture  and Open Problems} \label{5}
\bigskip

In spite of Theorem \ref{counterexample} and in view of Lemma \ref{Lemma of Alon}, we formulate a conjecture which sharpens the original Lemke-Kleitman Conjecture
for prime cyclic groups.

\begin{conj}
Let $G$ be a cyclic group of prime order and $S$ be a sequence over $G$ of length $|S| = |G|$.
Then $S$ has a subsequence $T$ with $\ind (T) = 1$ and length $|T| \in [1, \mathsf h(S)]$.
\end{conj}

\smallskip
\noindent
Let $G$ be a cyclic group of order $n \ge 2$. We denote by
\begin{itemize}
\item $\mathsf t (n)$ \ the smallest integer $\ell \in \mathbb N$ such that every sequence $S$ over $G$ of length $|S| \ge \ell$ has a subsequence $T$ with $\ind (T) = 1$.

\smallskip
\item $\mathsf T (n)$ \ the smallest integer $\ell  \in \mathbb N$ such that every squarefree sequence $S$ over $G$ of length $|S| \ge \ell$ has a subsequence $T$ with $\ind (T) = 1$.
\end{itemize}

\noindent
By Theorem \ref{counterexample}, it follows  that $\mathsf t(n)\geq n+\lfloor
\frac{n}{4} \rfloor -4$ for $n=4k+2\geq 22$.

\smallskip
\noindent
{\bf Open Problem.} {\sl Determine $\mathsf t(n)$ and $\mathsf T (n)$  for all  $n \ge 2$.}

\bigskip
\providecommand{\bysame}{\leavevmode\hbox to3em{\hrulefill}\thinspace}
\providecommand{\MR}{\relax\ifhmode\unskip\space\fi MR }
\providecommand{\MRhref}[2]{%
  \href{http://www.ams.org/mathscinet-getitem?mr=#1}{#2}
}
\providecommand{\href}[2]{#2}

\end{document}